\documentclass[12pt,reqno,a4paper,twoside]{article}
\usepackage{amsmath,amsthm,amstext,amscd,amssymb,euscript}
\usepackage{epsf}
\def\fee{\mathcal{F}}

\renewcommand{\phi}{\varphi}

\renewcommand{\emptyset}{\varnothing}

\def\limn{\lim_{n\to\infty}}

\def\1{ {\mathit{1} \!\!\>\!\! I} }
\newcommand\s{{\mathbf s}}

\makeatletter
\@addtoreset{equation}{section}
\makeatother

\newtheorem{ittheorem}{Theorem}
\newtheorem{itlemma}{Lemma}
\newtheorem{itproposition}{Proposition}
\newtheorem{itdefinition}{Definition}
\newtheorem{itremark}{Remark}

\newenvironment{theorem}{\addtocounter{equation}{1}
\begin{ittheorem}}{\end{ittheorem}}

\newenvironment{lemma}{\addtocounter{equation}{1}
\begin{itlemma}}{\end{itlemma}}

\newenvironment{proposition}{\addtocounter{equation}{1}
\begin{itproposition}}{\end{itproposition}}

\newenvironment{definition}{\addtocounter{equation}{1}
\begin{itdefinition}}{\end{itdefinition}}

\newenvironment{remark}{\addtocounter{equation}{1}
\begin{itremark}}{\end{itremark}}

\newcommand{\beq}{\begin{eqnarray}}
\newcommand{\eeq}{\end{eqnarray}}

\newcommand{\be}{\begin{equation}}
\newcommand{\ee}{\end{equation}}

\newcommand{\bl}{\begin{lemma}}
\newcommand{\el}{\end{lemma}}

\newcommand{\br}{\begin{remark}}
\newcommand{\er}{\end{remark}}

\newcommand{\bt}{\begin{theorem}}
\newcommand{\et}{\end{theorem}}

\newcommand{\bd}{\begin{definition}}
\newcommand{\ed}{\end{definition}}

\newcommand{\bp}{\begin{proposition}}
\newcommand{\ep}{\end{proposition}}

\newcommand{\bc}{\begin{corollary}}
\newcommand{\ec}{\end{corollary}}

\newcommand{\bpr}{\begin{proof}}
\newcommand{\epr}{\end{proof}}

\newcommand{\bi}{\begin{itemize}}
\newcommand{\ei}{\end{itemize}}

\newcommand{\ben}{\begin{enumerate}}
\newcommand{\een}{\end{enumerate}}

\newcommand{\R}{\mathbb R}
\newcommand{\N}{\mathbb N}

\newcommand{\E}{\mathbb E}

\newcommand{\pee}{\ensuremath{\mathbb{P}}}

\newcommand{\fe}{\ensuremath{\mathcal{F}}}

\newcommand{\si}{\ensuremath{\sigma}}

\newcommand{\qak}{\mathbb{Q}}
\newcommand{\de}{\delta}
\newcommand{\pot}{\pee^{\scriptscriptstyle{[0,t]}}}
\newcommand{\pott}{\peet^{\scriptscriptstyle{[0,t]}}}

\def\now{
\ifnum\time<60
          12:\ifnum\time<10 0\fi\number\time am
          \else
            \ifnum\time>719\chardef\a=`p\else\chardef\a=`a\fi
          \hour=\time
          \minute=\time
          \divide\hour by 60 
          \ifnum\hour>12\advance\hour by -12\advance\minute by-720 \fi
          \number\hour:%
          \multiply\hour by 60 
          \advance\minute by -\hour
          \ifnum\minute<10 0\fi\number\minute\a m\fi}
\newcount\hour
\newcount\minute
\numberwithin{equation}{section}         

\theoremstyle{remark}

\newcommand{\Lt}{\tilde{L}}
\newcommand{\ct}{\tilde{c}}
\newcommand{\peet}{\tilde{\pee}}
\newcommand{\mut}{\tilde{\mu}}
\newcommand{\pt}{\tilde{p}}

\begin{document}
\title{{\bf Relative entropy and waiting times for continuous-time
    Markov processes
}}
\author{
J.-R. Chazottes\footnote{CPhT, CNRS-Ecole polytechnique, 91128 Palaiseau
Cedex, France, and CMM, UMI CNRS 2807, Universidad de Chile,
Av. Blanco Encalada 2120, Santiago, Chile,
jeanrene@cpht.polytechnique.fr}\\
C. Giardina\footnote{Eurandom, Postbus 513,
5600 MB Eindhoven, The Netherlands, giardina@eurandom.tue.nl}\\
F.\ Redig\footnote{Mathematisch Instituut Universiteit Leiden,
Niels Bohrweg 1, 2333 CA Leiden, The Netherlands, redig@math.leidenuniv.nl}
}

\maketitle

\begin{abstract}
For discrete-time stochastic processes, there is a close connection
between return/waiting times and entropy. Such a connection cannot
be straightforwardly extended to the continuous-time
setting. Contrarily to the discrete-time case one does need a
reference measure and so the natural object is relative entropy rather
than entropy. In this paper we elaborate on this in the case of
continuous-time Markov processes with finite state space. A reference
measure of special interest is the one associated to the time-reversed
process. In that case relative entropy is interpreted as the entropy
production rate. The main results of this paper are: almost-sure
convergence to relative entropy of suitable waiting-times and their
fluctuation properties (central limit theorem and large deviation principle).

\end{abstract}

\section{Introduction}

Many limit theorems in the theory of stochastic processes have a version for discrete-time
as well as for continuous-time processes. The ergodic theory of Markov chains
e.g.\ is more or less identical in discrete and in continuous time. The same
holds for the Ergodic Theorem, martingale convergence theorems, central limit
theorems and large deviations for additive functionals, etc.
Usually, one obtains the same results with some additional effort in the continuous-time setting,
where e.g.\ extra measurability issues can show up.

For discrete-time ergodic processes, there is a remarkable theorem connecting
recurrence times and entropy \cite{shields}. In words, it states that
the logarithm of the first time the process repeats its first $n$ symbols typically behaves
like $n$ times the entropy of the process. This provides a way to sample entropy observing
a single, typical trajectory of the process.
This result seems a natural candidate to transport to a continuous-time setting.
The relation between entropy and return times is sufficiently intuitive
so that one would not expect major obstacles on the road toward such a result for
continuous-time ergodic processes.
There is however one serious problem. On the path
space of continuous-time processes (on a finite state space, say), there is no natural flat measure.
In the discrete-time setting one cannot distinguish between entropy of a process
and relative entropy between the process and the uniform measure on
trajectories. These only differ by a constant and a minus sign.
As we shall see below, this difference between relative entropy and entropy does play an important role
in turning to continuous-time processes.
Therefore, there will be no continuous-time analogue of the relation between return times
and ``entropy". In fact, the logarithm of return times turns out to have no suitable
way of being normalized, even for very simple processes in continuous time such as
Markov chains. To circumvent this drawback, we propose here is to consider {\em differences of
the logarithm of suitable waiting times} and relate them to {\em relative entropy}.

Another aspect of our approach is to first discretize time, and show that the relation
between waiting times and relative entropy persists in the limit of vanishing discrete time-step.
From the physical point of view, the time-step of the discretization
is the acquisition frequency of the device one uses to sample the process.
We can also think of numerical simulations for which the discretization of time is unavoidable.
Of course, the natural issue is to verify if the results obtained with the discretized process give the correct
ones for the original process, after a suitable rescaling and by letting the time-step go to zero.
This will be done in the present context.

In this paper, we will restrict ourselves to continuous-time Markov chains with finite-state space
for the sake of simplicity and also because the aforementioned problem already appears in this
special, yet fundamental, setting. The main body of this paper is: a law of large numbers
for the difference of the logarithm of certain waiting times giving a suitable relative entropy,
a large deviation result and a central limit theorem. One possible application is the estimation of
the relative entropy density between the forward and the backward
process which is physically interpreted as the {\em mean entropy
production}, and which is strictly
positive if and only if the process is reversible (i.e., in ``detailed balance'', or
``equilibrium'').

Our paper is organized as follows. In Section \ref{naive}, we show why the naive generalization
of the Ornstein-Weiss theorem fails. Section \ref{main} contains the main results about law of
large numbers, large deviations and central limit theorem for the logarithm of ratios of waiting times.
In the final section we consider the problem of ``shadowing" a given continuous-time trajectory drawn
from an ergodic distribution on path space.

\section{Naive approach}\label{naive}

In this section we start with an informal discussion motivating the quantities which we
will consider in what follows.
Let $\{X_t, t\geq 0\}$ be a continuous-time Markov chain with state space $A$,
with stationary measure $\mu$, and with generator
$$
Lf(x) = \sum_{y\in A} c(x) p(x,y) (f(y)-f(x))
$$
where $p(x,y)$ is a transition probability of a discrete-time irreducible Markov chain
on $A$, with $p(x,x)=0$, and where the escape rates $c(x)$ are strictly positive.
Given a time-step $\delta$, we can discretize the Markov chain to
obtain its ``$\delta$-discretization'' $\{X_{i\delta}, i=0,1,2,\ldots\}$.
Next we define the first time the $\delta$-discretized process repeats
its first $n$ symbols via the random variable
\be\label{hitler}
R^\delta_n (X) : =
\inf\{ k\geq 1: (X_{0},\ldots,X_{(n-1)\delta})= (X_{k\delta},\ldots, X_{(k+n-1)\delta})\}\,.
\ee
The analogue of the Ornstein-Weiss theorem \cite{shields} for this continuous-time process would be
a limit theorem for a suitably normalized version of  $\log R^\delta_n$ for $n\to\infty, \delta\to 0$.
However, for $\delta>0$ fixed, the ergodicity of the
$\delta$-discretization $\{ X_0, X_{\delta}, X_{2\delta}, \ldots,
X_{n\delta},\ldots\}$, Ornstein-Weiss and Shannon-McMillan-Breiman
theorems \cite{shields} yield
$$
\frac1n \log\left[R^\delta_n (X) \pee(X_1^n)\right]=
o(1)\quad\textup{eventually a.s. as }\, n\to\infty\,.
$$
Using the fact that $X$ is an ergodic Markov chain we obtain
\beq\label{ow}
-\frac1n \log
R^\delta_n (X) &=&
\frac1n\log\mu (X_0)+ \frac1n \left( \sum_{i=1}^{n-1} \log p^X_\delta (X_i, X_{i+1}) \right) +o(1)
\nonumber\\
& = &\E \left[ \log p^X_\delta (X_{0}, X_{1}) \right] + o(1)
\quad\textup{eventually a.s. as }\ n\to\infty
\eeq
where $\E$ denotes expectation in the Markov chain started from its
stationary distribution and where $p^X_\delta$ denotes the transition probability of the
$\delta$-discretized Markov chain $\{X_{i\delta}, i=0,1,2,\ldots\}$, i.e.,
$$
p^X_\delta(x,y) = (e^{\delta L})_{xy}
=
\1_{xy}(1-\delta c(x)) + \delta c(x) p(x,y)  + \mathcal{O} (\delta^2)\,.
$$
Therefore,
\beq
&& -\lim_{n\to\infty}\frac1n \log R^\delta_n (X)=
\nonumber\\
& &\sum_{x\in A} \mu (x) (1-\delta c(x)) \log (1-\delta c(x)) +
\sum_{x,y\in A}\mu (x) \delta c(x) p(x,y) \log(\delta c(x) p(x,y))+ \mathcal{O}(\delta^2)\,.
\nonumber
\eeq

In this expression we see that the first term is of order $\delta$ whereas the second one is of order
$\delta\log\delta$. Therefore, this expression does not seem to have a natural way to be normalized.
This is a typical phenomenon for continuous-time processes: we need a suitable reference
process in order to define ``entropy" as ``relative entropy" with respect to this reference process.
Indeed, as we will see in the next sections, by considering {\em differences of waiting times}
one is able to cancel the $\delta\log\delta$ term in order to obtain expression that
makes sense in the limit $\delta\downarrow 0$.

\section{Main results: waiting times and relative entropy}\label{main}

We consider continuous-time Markov chains with a finite state-space $A$. We will
always work with irreducible Markov chains with a unique stationary distribution.
The process is denoted by $\{ X_t:t\geq 0\}$. The associated measure on path
space starting from $X_0=x$ is denoted by $\pee_x$ and by $\pee$ we
denote the path space measure of the process started from its unique stationary distribution.
For $t\geq 0$, $\fe_t$ denotes the sigma-field generated by $X_s,s\leq t$, and $\pee^{[0,t]}$ denotes
the measure $\pee$ restricted to $\fe_t$.

\subsection{Relative entropy: comparing two Markov chains}\label{RE}

Consider two continuous-time Markov chains, one denoted by
$\{ X_t: t\geq 0\}$ with generator
\be\label{gen1}
Lf(x) = \sum_{y\in A} c(x) p(x,y) (f(y)-f(x))
\ee
and the other denoted by $\{ Y_t: t\geq 0 \}$ with generator
$$
\tilde{L} f(x) = \sum_{y\in A} \tilde{c}(x) p(x,y) (f(y)-f(x))
$$
where $p(x,y) $ is the Markov transition function of an
irreducible discrete-time Markov chain.
We further assume that $p(x,x)=0$, and $c(x)>0$ for all $x\in A$.
We suppose that $X_0$, resp.\ $Y_0$, is distributed according to the unique
stationary measure $\mu$, resp $\tilde{\mu}$
so that both processes are stationary and ergodic.

\br
The fact that the Markov transition function $p(x,y)$ is the same
for both processes is only for the sake of simplicity. All our results
can be reformulated in the case that the Markov transition functions
would be different.
\er

We recall Girsanov's formula \cite{GS}:
\be\label{grizzly}
\frac{d\pot}{d\pott} (\omega)
=\frac{\mu(\omega_0)}{\mut_(\omega_0)}
\exp\left( \int_0^t \log\frac{c(\omega_s)}{\ct (\omega_s)} \ dN_s(\omega) -
  \int_0^t (c(\omega_s) -\ct(\omega_s))\ ds\right)
\ee
where $N_s(\omega)$ is the number of jumps of the path $\omega$ up to time $s$.
The relative entropy of $\pee$ w.r.t. $\peet$ up to time $t$ is
defined as
\be\label{relen}
s_t (\pee|\peet) = \int d\pee (\omega) \log \left(\frac{d\pot}{d\pott} (\omega)\right)\,.
\ee

Using \eqref{grizzly} and stationarity, we obtain
\beq\label{relen2}
\lim_{t\to\infty}\frac{s_t (\pee|\peet)}{t} & = & \sum_{x\in A} \mu (x)
  c(x) \log\frac{c(x)}{\ct (x)} - \sum_{x\in A} \mu (x) (c(x)-\ct
  (x))\\
&=: & \s(\pee|\peet)
\nonumber
\eeq
where $\s(\pee|\peet)$ is the {\em relative entropy (per unit time)} of
$\pee$ with respect to $\peet$. We refer to \cite{DZ}, \cite{varadhan} for more details on relative
entropy for continuous-time Markov chains.
Notice also that, by ergodicity,
$$
\lim_{t\to\infty}\frac1t \log \frac{d\pot}{d\pott}
(\omega) =\s(\pee|\peet)\quad\pee-\textup{a.s.}\,.
$$
In the case $\{ Y_t:t\geq 0\}$ is Markov chain with generator
$$
\tilde{L} f(x) = \sum_{y\in A} \tilde{c}(x) \pt(x,y) (f(y)-f(x))
$$
\eqref{relen2} generalizes to
\be
s(\pee|\peet)
=
\sum_{x,y\in A} \mu (x)
  c(x) p(x,y) \log\frac{c(x)p(x,y)}{\ct (x)\pt(x,y)} - \sum_{x\in A} \mu (x) (c(x)-\ct
  (x))
\ee
A important particular case is met when $\{ Y_t:t\geq 0\}$ is the time-reversed
process of $\{ X_t:t\geq 0\}$, i.e.,
\[
(Y_t)_{0\leq t\leq T} = (X_{T-t})_{0\leq t\leq T}\ \mbox{in\  distribution}\,.
\]
This is a Markov chain with transition rates
\be\label{rever}
\ct (x,y) = c(x) \frac{c(y) p(y,x) \mu (y) }{c(x)\mu (x)} \,.
\ee
In that particular situation, the  random variable
\be\label{entpro}
S_T(\omega) = \log \frac{d\peet^{[0,T]}}{d\pee^{[0,T]}}
\ee
has the interpretation of ``entropy production", and the relative entropy density
$s(\pee|\peet)$ has the interpretation of ``mean entropy production per unit time".
see e.g. \cite{JQQ,maes}.

\subsection{Law of large numbers}\label{LLN}

For $\delta>0$, we define the discrete-time Markov chain $X^\delta:=
\{ X_0, X_\delta, X_{2\delta}, \ldots\}$. This Markov chain
has transition probabilities
\beq\label{pxdel}
p^X_\delta(x,y) &=& (e^{\delta L})_{xy}
\nonumber\\
&=&
\1_{xy}(1- \delta c(x)) + \delta c(x) p(x,y) + \mathcal{O}(\delta^2)
\eeq
where $\1$ is the identity matrix.
Similarly we define another Markov chain $Y^\delta$ with transition probabilities
\beq\label{py}
p^Y_\delta(x,y) &=& (e^{\delta \Lt})_{xy}
\nonumber\\
&=&
\1_{xy} (1- \delta \ct(x))+ \delta \ct(x)p(x,y) + \mathcal{O}(\delta^2)\,.
\eeq
The path-space measure (on $A^\N$) of $X^\delta$, resp.\ $Y^\delta$,
is denoted by $\pee^\delta$, resp.\ $\peet^\delta$. From now on,
we will write $\pee^\delta(X_1^n)$ instead of
$\pee^\delta(X_\delta,X_{2\delta}\ldots,X_{n\delta})$ to alleviate notations.

We define waiting times, which are random variables
defined on $A^\N \times A^\N$, by setting
\be\label{wxy}
W_{n}^{\delta} (X|Y) = \inf \{ k\geq 1: (X^\delta_1,\ldots,X^\delta_n) = (Y^\delta_{k+1},\ldots,Y^\delta_{k+n})\}
\ee
where we make the convention $\inf\emptyset =\infty$.
In words, this is the first time that in a realization of the process
$Y^\delta$ that one observes the first
$n$ symbols of a realization of the process $X^\delta$.
Similarly, if ${X'}^\delta$ is an independent copy of the process $X^\delta$, we define
\be\label{wxx}
W_{n}^{\delta} (X|X') = \inf \{ k\geq 1: (X^\delta_1,\ldots,X^\delta_n) = (X'^\delta_{k+1},\ldots,X'^\delta_{k+n})\}\,.
\ee
We then have the following law of large numbers.

\bt\label{lln}
($\pee\otimes\peet\otimes\pee$)-almost surely:
\beq\label{drac}
&&
\lim_{\delta\to 0}\lim_{n\to\infty}
\frac1{n\delta} \log\frac{W_{n}^{\delta} (X|Y)}{W_{n}^{\delta} (X|X')}
=\s(\pee|\peet)\,.
\eeq
\et

Before proving this theorem, we state a theorem about the exponential
approximation for the hitting-time law, which will be the crucial ingredient
throughout this paper.
For a $n$-block $x_1^n:=x_1,\ldots,x_n\in A^n$ and a discrete-time
trajectory $\omega\in A^\N$, we define the hitting time
\be\label{hit}
T^\delta_{x_1^n} (\omega) = \inf\{ k\geq 1: X_{(k+1)\delta}=\omega_1,\ldots,X_{(n+k+1)\delta}=\omega_{n+1}\}\,.
\ee
We then have the following result, see \cite{abadi}.

\bt\label{explaw}
For all $\delta >0$, there exist $\eta_1,\eta_2, C, c,\beta,\kappa\in ]0,\infty[$ such that
for all $n\in\N$ and for all $x_1^n\in A^n$, there exists $\eta= \eta(x_1^n)$, with
$0<\eta_1\leq \eta\leq\eta_2<\infty$ such that for all $t>0$
\beq\label{expform}
\Big |\pee \left(T^\delta_{x_1^n} (\omega)>\frac{t}{\pee^\delta(X_1^n=x_1^n)}\right) -e^{-\eta t}\Big|
&\leq &
Ce^{-ct} \left(\pee^\delta (X_1^n=x_1^n)\right)^\kappa
\nonumber\\
&\leq &
Ce^{-ct}e^{-\beta n}\,.
\eeq
The same theorem holds with $\pee$ replaced by $\peet$.
\et

The constants appearing in Theorem \ref{explaw} (except $C$) depend
on $\delta$, and more precisely we have $\beta=\beta(\delta)\to 0$,
$\eta_1 =\eta_1 (\delta) \to 0$ as $\delta\to 0$.

This is important in applications, since one wants to choose a certain
discretization $\delta$ and then a corresponding ``word-length" $n(\delta)$ for the waiting times, or vice-versa.
From Theorem \ref{explaw} we derive (see \cite{acrv}):

\bp\label{sandwich}
For all $\delta>0$, there exist $\kappa_1,\kappa_2 >0$ such that
$$
-\kappa_1\log n\leq \log \left(W_{n}^{\delta} (X|Y) \peet^\delta (X_1^n)\right)\leq
\log(\log n^{\kappa_2})\quad\pee\otimes\peet \ \textup{eventually a.s.}
$$
and
$$
-\kappa_1\log n\leq \log \left(W_{n}^{\delta} (X|X') \pee^\delta ({X}_1^{n})\right)\leq
\log(\log n^{\kappa_2})\quad\pee \otimes\pee \ \textup{eventually a.s.}
$$
\ep

With these ingredients we can now give the proof of Theorem \ref{lln}.

\bigskip

{\em Proof of Theorem \ref{lln}}.
From Proposition \ref{sandwich} it follows
that, for all $\delta>0$, $\pee\otimes\peet \otimes\pee$ almost surely
\be
\lim_{n\to\infty}\frac{1}{n}\left(\log W_{n}^{\delta} (X|Y) - \log W_{n}^{\delta} (X|X')
+\sum_{i=0}^{n-1} \log p^Y_\delta (X_i,X_{i+1})
-
\sum_{i=0}^{n-1} \log p^X_\delta (X_i,X_{i+1})\right)=0\,.
\ee
By ergodicity of the continuous-time Markov chain $\{ X_t:t\geq 0\}$, the discrete
Markov chains $X^\delta, Y^\delta$ are also ergodic and therefore we obtain
\be
\lim_{n\to\infty}\frac1{n}\left(\log W_{n}^{\delta} (X|Y) - \log W_{n}^{\delta} (X|X')
+\sum_{x,y\in A} \mu(x) p^X_{\delta} (x,y) \log \left(\frac{p^Y_\delta (x,y)}{p^X_\delta (x,y)}\right)
\right)=0\,.
\ee
Using
\eqref{pxdel}, \eqref{py} and $p(x,x)=0$, this gives
\beq
&&\lim_{n\to\infty}\frac1n\left(\log W_{n}^{\delta} (X|Y) - \log W_{n}^{\delta} (X|X')\right)
\nonumber\\
&=&
-\sum_{x\in A} \mu(x) (1-\delta c(x)) \log (1-\delta \ct (x))
-\sum_{x,y\in A} \mu (x) \delta c(x) p(x,y) \log (\delta \ct (x) p(x,y))
\nonumber\\
&+&
\sum_{x\in A} \mu(x) (1-\delta c(x)) \log (1-\delta c (x))
+\sum_{x,y\in A} \mu (x) \delta c(x) p(x,y) \log (\delta c (x) p(x,y))
+ \mathcal{O}(\delta^2)
\nonumber\\
&=&\delta
\left(\sum_{x,y\in A} \mu(x) c(x) p(x,y) \log \frac{c(x)}{\ct (x)}
  +\sum_{x\in A}\mu (x) (\ct (x)-c(x))\right) + \mathcal{O}(\delta^2)
\nonumber\\
&=&\delta\ \s(\pee|\peet)+ \mathcal{O}(\delta^2)\,.
\eeq
Combining this with \eqref{relen2} concludes the proof of Theorem \ref{lln}. \hfill$\square$

Let us now specify the dependence on $\delta$ of the various constants appearing
in Theorem \ref{explaw}.
For the lower bound on the parameter we have (see \cite{abadi}, section 5)
\be\label{blub}
\eta_1 (\delta) \geq \frac{1}{C'+K}
\ee
where $C'$ is a positive number independent of $\delta$ and
$$
K= 2\sum_{l=1}^\infty \alpha (l) + \sum_{k=1}^{n/2} \sup_{\{x_1^{(n-k)}\}} \pee^\delta (X_1^{(n-k)}=x_1^{(n-k)})\,.
$$
Here $\alpha(l)$ denotes the classical $\alpha$-mixing coefficient:
\[
\alpha(l) = \sup_{j\geq 1} \sup_{S_1\in \fe_0^{j-1},
  S_2\in\fe_{j+l}^\infty}
\left(\pee^\delta (S_1\cap S_2) - \pee^\delta (S_1)\pee^\delta (S_2)\right)
\]
where $\fe_{m}^{n}$ is the Borel sigma-field on $A^\N$ generated by
$X_m^n$ ($0\leq m\leq n\leq \infty$).
By the assumption of ergodicity of the continuous Markov chain, the generator $L$ (resp. $\Lt$) has an eigenvalue
$0$, the largest real part of the other eigenvalues is strictly
negative and denoted by $-\lambda_1 <0$, and one has
\be\label{mix}
\alpha(l) \leq \exp (-\lambda_1\delta l)\, .
\ee
Using \eqref{pxdel} there exists $\lambda_2>0$ such that
$$
\pee^\delta (X_1^{(n-k)}=x_1^{(n-k)})\leq \exp(-\lambda_2\delta n/2)
$$
for $k=1,\ldots,n/2$.
Therefore, there exists $\hat{c} >0$ such that
\be\label{la1}
\eta_1 (\delta) > \hat{c}\delta\,.
\ee
Similarly, from the proof of Theorem 2.1 in \cite{acrv}  one obtains easily the dependence on $\delta$
of the  constants appearing in the error term of
\eqref{expform}.
\be\label{error}
c= c(\delta) > \gamma_1\delta, \beta=\beta(\delta) > \gamma_2 \delta
\ee
for some $\gamma_1,\gamma_2 >0$.

In applications, e.g., the estimation of the relative entropy from a sample path,
one would like to choose the word-length $n$ and the discretization
$\delta=\delta_n$ together.
This possibility is precisely provided by the estimates \eqref{la1}
and \eqref{error}, as the following analogue of Proposition \ref{sandwich} shows.
\bp\label{sandwichn}
Let $\delta_n\to 0$ as $n\to \infty$, then
there exists $\kappa_1,\kappa_2 >0$
\be\label{bam}
-\kappa_1\frac{\log n}{\delta_n}\leq \log \left(W_{n}^{\delta_n} (X|Y) \peet^{\delta_n }(X_1^n)\right)\leq
\log\left(\frac{\kappa_2\log n}{\delta_n}\right)\quad\pee\otimes\peet\ \textup{eventually a.s.}
\ee
and
$$
-\kappa_1\frac{\log n}{\delta_n}\leq \log \left(W_{n}^{\delta_n} (X|X') \pee^{\delta_n} ({X}_1^n)\right)\leq
\log\left(\frac{\kappa_2\log n}{\delta_n}\right)\quad\pee\otimes\pee\ \textup{eventually a.s.}\,.
$$
\ep
\bpr
The proof is analogous to the proof of Theorem 2.4 in \cite{acrv}.
For the sake of completeness, we prove the upper bound \eqref{bam}.
We can assume that $\delta_n\leq 1$.
By the exponential approximation \eqref{expform} we have, for all
$t>0$, $n\geq 1$, the estimates
\beq\label{boem}
\pee^\delta \otimes\peet^\delta\left(\log \left(W_{n}^{\delta} (X|Y) \peet^\delta (X_1^n)\right)\geq \log t\right)
&\leq &
e^{-\eta (\delta_n) t} + C e^{-\beta (\delta_n)n }e^{-c(\delta_n)t}
\nonumber\\
\leq
e^{-\eta_1 \delta_n t} + C e^{-\gamma_1 \delta_n n}e^{-\gamma_2 \delta_n t}\,.
\eeq
Choosing $t=t_n = \frac{\kappa_2\log n}{\delta_n}$, with $\kappa_2>0$
large enough makes the rhs of \eqref{boem} summable
and hence a Borel-Cantelli argument gives the upper bound.
\epr
Of course, whether this proposition is still useful, i.e., whether it still gives the law of large numbers
with $\delta=\delta_n$ depends on the behavior of ergodic sums
\[
\sum_{i=1}^n f(X_i)
\]
under the measure $\pee^{\delta_n}$, i.e., the behavior of
\[
\sum_{i=1}^n f(X_{i\delta_n})
\]
under $\pee$.
This is made precise in the following theorem:
\bt\label{lln1} Suppose that $\delta_n\to 0$ as $n\to\infty$ such that
$\frac{\log n}{n\delta^2_n}\to 0$ then in
($\pee\otimes\peet\otimes\pee$) probability:
\beq\label{drac1}
&&
\lim_{n\to\infty}
\frac1{n\delta_n} \log\frac{W_{n}^{\delta_n} (X|Y)}{W_{n}^{\delta_n} (X|X')}
=\s(\pee|\peet)\,.
\eeq
\et
\bpr
By proposition \ref{sandwichn} we can write
\beq
\label{qaws}
&&\log W^{\delta_n}_n (X|Y)-\log W^{\delta_n}_n (X|X')
\nonumber\\
&=&
\sum_{i=1}^n \1_{X_i=X_{i+1}}\log\frac{1-\delta_n c(X_i)}{1-\delta_n \tilde{c} (X_i)}
+
\sum_{i=1}^n \1_{X_i\not=X_{i+1}}\log\frac{c(X_i)}{\tilde{c} (X_i)}
+ \mathcal{O}(\log n/\delta_n)\,.
\eeq
The sum on the right hand site of \eqref{qaws} is of the form
\be\label{qwer}
\sum_{i=1}^n F_{\delta_n} (X_{i\delta_n}, X_{(i+1)\delta_n})
\ee
with
\be
\E( F_{\delta_n}-\E(F_{\delta_n}))^2 \leq C\delta_n
\ee
where $C>0$ is some constant. Now, using ergodicity of the
continuous-time Markov chain $\{X_t, t\geq 0\}$, we have the estimate
\be
\E\Big((F_{\delta_n} (X_{i\delta_n}, X_{(i+1)\delta_n})-\E(F_{\delta_n}))
(F_{\delta_n} (X_{j\delta_n}, X_{(j+1)\delta_n})-\E(F_{\delta_n}))
\Big)
\leq
\|F_{\delta_n }\|_2^2 \ e^{-\delta_n\lambda_1 |i-j|}
\ee
with $\lambda_1>0$ independent of $n$.

Combining these estimates gives
\be
Var\left(\sum_{i=1}^n F_{\delta_n} (X_{i\delta_n}, X_{(i+1)\delta_n})\right)
\leq Cn\delta_n + \sum_{i=1}^n \sum_{j\in\{1,\ldots,n\}\setminus\{i\}} \delta_n e^{-\delta_n\lambda_1 |i-j|}
\leq C n\delta_n + C' \delta_n\frac{n}{\delta_n}
\ee
where $C'>0$ is some constant.
Therefore,
\be\label{qasd}
\frac{1}{n^2\delta_n^2}
Var\left(\sum_{i=1}^n F_{\delta_n} (X_{i\delta_n}, X_{(i+1)\delta_n})\right)
= \mathcal{O}(1/n\delta_n^2)\,.
\ee
Combining
\eqref{qaws} and \eqref{qasd} with the assumption $\frac{\log n}{n\delta^2_n}\to 0$ concludes the proof.
\epr

\subsection{Large deviations}

In this subsection, we study the large deviations of
\[
\frac{1}{n}\log\left(\frac{W_{n}^{\delta} (X|Y)}{W_{n}^{\delta} (X|X')}\right)\, \cdot
\]
More precisely, we compute the large deviation generating function $\fe^\delta (p)$ in the
limit $\delta\to 0$ and show that it coincides with the large deviation
generating function for the Radon-Nikodym derivatives $d\pee^{[0,t]}/d\peet^{[0,t]}$.
As in the case of waiting times for discrete-time processes,
see e.g.\ \cite{cr}, the scaled-cumulant generating function
is only finite in the interval $(-1,1)$.

For the sake of convenience we introduce the function
$$
\mathcal{E}(p):= \lim_{\delta\to 0} \mathcal{E}^\delta(p)=
\lim_{t\to\infty}\frac1t
\log \E_{\pee}
\left( \frac{d\pot}{d\pott}\right)^p =
$$
\beq\label{genfun}
\lim_{t\to\infty}\frac1t
\log \E_{\pee}
\left(\exp\left(
p\left(
\int_0^t \log\frac{c(\omega_s)}{\ct (\omega_s)} \ dN_s(\omega) -
  \int_0^t (c(\omega_s) -\ct(\omega_s))\ ds\right)\right)\right)\,.
\eeq
By standard large deviation theory for continuous-time Markov chains
(see e.g. \cite{stroock}) this function exists and is the scaled-cumulant generating function for
the large deviations of
\[
\int_0^t \log\frac{c(\omega_s)}{\ct (\omega_s)} \ dN_s(\omega) -
  \int_0^t (c(\omega_s) -\ct(\omega_s))\ ds
\]
as $t\to\infty$.

We can now formulate the following large deviation theorem.

\bt\label{largedev}
For all $p\in \R$ and $\delta>0$ the function
\be\label{karamel}
\fee^\delta (p):= \lim_{n\to\infty}
\frac1{n\delta} \log
\E_{\pee^{\delta}\otimes\peet^{\delta}\otimes\pee^{\delta}}\left( \frac{W_{n}^{\delta} (X|Y)}{W_{n}^{\delta} (X|X')}\right)^p
\ee
exists, is finite in $p\in (-1,1)$ whereas
$$
\fee^\delta(p)=\infty\quad\textup{for}\; |p|\geq 1\, .
$$

Moreover, as $\delta\to 0$, we have, for all $p\in (-1,1)$:
$$
\fee (p):=\lim_{\delta\to 0} \fee^\delta (p)=\mathcal{E}(p)\,.
$$
\et

The following notion of logarithmic equivalence will be convenient later on.
\bd\label{bombom}
Two non-negative sequences $a_n$, $b_n$ are called {\em logarithmically equivalent} (notation
$a_n\simeq b_n$)
if
\[
\lim_{n\to\infty} \frac1n (\log a_n -\log b_n) =0\,.
\]
\ed

\bpr
To prove  Theorem \ref{largedev}, we start with the following lemma.
\bl
\ben
\item
For all $\delta>0$ and for $|p|<1$,
\be\label{pape}
\E_{\pee^{\delta}\otimes\peet^{\delta}\otimes\pee^{\delta}}\left( \frac{W_{n}^{\delta} (X|Y)}{W_{n}^{\delta} (X|X')}\right)^p
\simeq
\E_{\pee^\delta}
\exp\left(p\sum_{i=0}^{n-1} \log\left(\frac{p^X_\delta (X_i,X_{i+1})}{p^Y_\delta (X_i,X_{i+1})}\right)\right)\,.
\ee
\item
For $|p|>1$,
\be\label{divpo}
\lim_{n\to\infty}\frac1n\log\E_{\pee^{\delta}\otimes\peet^{\delta}\otimes\pee^{\delta}}\left(
\frac{W_{n}^{\delta} (X|Y)}{W_{n}^{\delta} (X|X')}\right)^p
=\infty\,.
\ee
\een
\el
\bpr
The proof is similar to that of Theorem 3 in \cite{cr}.
\beq
&  &\E_{\pee^{\delta}\otimes\peet^{\delta}\otimes\pee^{\delta}} \left(
\frac{W_{n}^{\delta} (X|Y)}{W_{n}^{\delta} (X|X')}\right)^p
\nonumber\\
&=&
\sum_{x_1,\ldots,x_n}
\pee^\delta (X_1^n=x_1^n)
\left(
\frac{\pee^\delta({X'}_1^n=x_1^n)}{\peet^\delta(Y_1^n=x_1^n)}\right)^p
\nonumber\\
& \times &
\E_{\peet^\delta\otimes\pee^\delta}
\left(\frac{T_{x_1^n} (Y^\delta)\peet^\delta(Y_1^n=x_1^n)}{T_{x_1^n} (X'^\delta)\pee^\delta({X'}_1^n=x_1^n)}
\right)^p
\nonumber\\
&=&
\sum_{x_1,\ldots,x_n}
\pee^\delta (X_1^n=x_1^n)^{1+p}
\, \peet^\delta(Y_1^n=x_1^n)^{-p}
\, \E_{\peet^\delta\otimes\pee^\delta}
\left(\frac{\xi_n}{\zeta_n}\right)^p
\nonumber
\eeq
where
\[
\xi_n= T_{x_1^n} (Y^\delta)\, \peet^\delta(Y_1^n=x_1^n)
\]
and
\[
\zeta_n=T_{x_1^n} (X'^\delta)\, \pee^\delta({X'}_1^n=x_1^n)\,.
\]
The random variables $\xi_n,\zeta_n$ have approximately an exponential distribution
(in the sense of Theorem \ref{explaw}) and are independent.
Using this, we can repeat the arguments of the proof
of Theorem 3 in \cite{cr} -which uses the exponential law with the error-bound given
by Theorem \ref{explaw}- to prove that for $p\in (-1,1)$
\[
0< C_1\leq \E_{\peet^\delta\otimes\pee^\delta}
\left(\frac{\xi_n}{\zeta_n}\right)^p \leq C_2 <\infty
\]
where $C_1,C_2$ do not depend on $n$, whereas for $|p|>1$,
\be\label{divex}
\E_{\peet^\delta\otimes\pee^\delta}
\left(\frac{\xi_n}{\zeta_n}\right)^p =\infty\,.
\ee

Therefore, with the notation of Definition \ref{bombom}, for $|p|<1$
\beq\label{wawawa}
&&\E_{\pee^{\delta}\otimes\peet^{\delta}\otimes\pee^{\delta}} \left( \frac{W_{n}^{\delta} (X|Y)}{W_{n}^{\delta} (X|X')}\right)^p
\nonumber\\
&\simeq &
\sum_{x_1,\ldots,x_n}
\pee^\delta (X_1=x_1,\ldots,X_n=x_n)^{1+p}
\, \peet^\delta(Y_1=x_1,\ldots,Y_n=x_n)^{-p}
\nonumber\\
&=&
\E_{\pee^\delta}
\exp\left(p\sum_{i=1}^n \log\left(\frac{p^X_\delta (X_i,X_{i+1})}{p^Y_\delta (X_i,X_{i+1})}\right)\right)\,\cdot
\eeq
and for $|p|>1$
we obtain \eqref{divpo} from \eqref{divex}.
\epr
This proves the existence of $\fe^\delta (p)$. Indeed, the limit
\be\label{wawa}
\fe^\delta (p)= \lim_{n\to\infty}\frac1{n}\log\E_{\pee^\delta}
\exp\left(p\sum_{i=1}^n \log\left(\frac{p^X_\delta (X_i,X_{i+1})}{p^Y_\delta (X_i,X_{i+1})}\right)\right)
\ee
exists by standard large deviation theory of (discrete-time, finite state space) Markov chains
(since $\delta>0$ is {\em fixed}).

In order to deal with the limit $\delta\to 0$ of $\fe^\delta (p)$,
we expand the expression in the rhs of \eqref{pape}, up to order $\delta^2$. This gives
\beq
&&\E_{\pee^\delta} \exp \left(p\sum_{i=1}^n \log\left(\frac{p^X_\delta (X_i,X_{i+1})}{p^Y_\delta (X_i,X_{i+1})}\right)\right)
\nonumber\\
=&&
e^{\mathcal{O}(n\delta^2)}\E_{\pee^\delta}\left(\exp\left(
p\sum_{i=1}^n \log\frac
{\1_{X_i,X_{i+1}} + \delta c(X_i) p(X_i,X_{i+1}) - \delta c(X_i)}{\1_{X_i,X_{i+1}} +
\delta \ct(X_i) p(X_i,X_{i+1}) - \delta \ct(X_i)}\right)\right)
\nonumber
\eeq
\beq
\!\!\!\!\!\!\!\!\!\!\!\!\!=e^{\mathcal{O}(n\delta^2)}\E_{\pee^\delta}\Big[\exp
\Big(
p\sum_{i=1}^n \delta \1(X_i=X_{i+1}) \!\!\!\!\!\!\!\!\!\!\!\!\!& \!\!\!\!\!\!\!\!\!\!\!\!\!& (\ct (X_i) -c(X_i))
\Big)
\nonumber \\
& + & p\sum_{i=1}^n \1(X_i\not= X_{i+1} )\log\frac{c(X_i)}{\ct(X_i)}\, \Big]
\nonumber
\eeq
\beq
\quad =e^{\mathcal{O}(n\delta^2)}\E_{\pee}\Big[\exp
\Big(
p\sum_{i=1}^n \delta \1(X_{i\delta}=X_{(i+1)\delta})\!\!\!\!\!\!\!\!\!\!\!\!\!& \!\!\!\!\!\!\!\!\!\!\!\!\!&
(\ct (X_{i\delta}) -c(X_{i\delta}))
\Big)
\nonumber \\
& + & p\sum_{i=1}^n \1(X_{i\delta}\not= X_{(i+1)\delta} )\log\frac{c(X_{i\delta})}{\ct(X_{i\delta})}\, \Big] \, \cdot
\nonumber
\eeq

Next we prove that for all $K\in \R$
\be\label{karate}
\small{
\log\E_{\pee^\delta}\left[
\frac{
\exp
\left(
K\sum_{i=1}^n \big(
\delta \1(X_{i\delta}=X_{(i+1)\delta}) (\ct (X_{i\delta}) -c(X_{i\delta}))
+ \1(X_{i\delta}\neq X_{(i+1)\delta} )\log\frac{c(X_{i\delta})}{\ct(X_{i\delta})}\big)\right)
}
{
\exp\left(K\int_0^{n\delta} (\ct(X_s)-c(X_s))ds + K \int_0^{n\delta}\log\frac{c(X_s)}{\ct(X_s)}dN_s\right)
}
\right]
}
\ee
$
=\mathcal{O}(n\delta^2)
$.

\bigskip

This implies the result of the theorem by a standard application of H\"{o}lder's inequality, see
e.g., \cite{DZ}.
We first consider the difference
$$
A(n,\delta):=\Big |\sum_{i=1}^n \1(X_{i\delta}\not= X_{(i+1)\delta} )\log\frac{c(X_{i\delta})}{\ct(X_{i\delta})}
-\int_0^{n\delta} \log\frac{c(X_s)}{\ct(X_s)}dN_s\Big|\, .
$$

If there does not exist an interval $[i\delta,(i+1)\delta[, i\in \{0,\ldots,n-1\}$ where at least two jumps of the Poisson
process $\{N_t, t\geq 0\}$ occur, then $A(n,\delta)=0$.
Indeed, if there is no jump in $[i\delta,(i+1)\delta[$, both
$ \1(X_{i\delta}\not= X_{(i+1)\delta} )\log\frac{c(X_{i\delta})}{\ct(X_{i\delta})}$
and
$\int_{i\delta}^{(i+1)\delta} \log\frac{c(X_s)}{\ct(X_s)}dN_s$
are zero and if there is precisely one jump, then they are equal.
Therefore, using the independent increment property of the Poisson process, and the strict positivity of the rates,
we have the bound
$$
A(n,\delta) \leq C\sum_{i=1}^n \1(\chi_i\geq 2)
$$
where the $\chi_i$´s, $i=1,\ldots,n$, form a collection of independent Poisson random variables with
parameter $\delta$, and $C$ is some positive constant.
This gives
\be\label{and}
\E_{\pee^\delta} e^{2K A(n,\delta)} = (\mathcal{O}(\delta^2) e^{2K} + \mathcal{O}(1))^n = \mathcal{O}(e^{n\delta^2})\,.
\ee
Next, we tackle
\be\label{bnd}
B(n,\delta) :=
\Big|\sum_{i=1}^n \delta \1(X_{i\delta}=X_{(i+1)\delta}) (\ct (X_{i\delta}) -c(X_{i\delta}))
-
\int_0^{n\delta} (\ct(X_s)-c(X_s))ds
\Big|\, .
\ee
If there is no jump in any of the intervals $[i\delta, (i+1)\delta[$, this term is zero.
Therefore is is bounded by
$$
B(n,\delta) \leq
C'\delta \sum_{i=1}^n \1(\chi_i\geq 1)
$$
where the $\chi_i$´s, $i=1,\ldots,n$, form once more a collection of independent Poisson random variables with
parameter $\delta$, and $C'$ is some positive constant.
This gives
$$
\E_\pee e^{2K B(n,\delta)} \leq ( \mathcal{O}(\delta e^{C''\delta}) + 1-\delta)^n
= \mathcal{O}(e^{n\delta^2})
$$
where $C''$ is some positive constant.
Hence, \eqref{karate} follows by combining \eqref{and} and \eqref{bnd} and using Cauchy-Schwarz inequality.
\epr
The following propoisition is a straightforward application of Theorem \ref{largedev}
and \cite{PS}.
\bp
For all $\delta >0$, $\fee^\delta$ is real-analytic and convex, and
the sequence $\{ \log W_{n}^{\de}(X|Y)-\log W_{n}^{\de}(X|X'): n\in \N\}$
satisfies the following large large deviation principle: Define the
open interval $(c_-,c_+)$, with
\[
c_{\pm}:=\lim_{p\to \pm 1}\frac{d\mathcal{E}^{\delta}}{dp}<0
\]
Then, for every interval $J$
such that $J\cap (c_-, c_+)\neq \emptyset$
$$
\limn \frac{1}{n}\log \pee^{\delta}\otimes\peet^{\delta}\otimes\pee^{\delta}\left\{
\frac{1}{n}\log\left(\frac{W_{n}^{\de}(X|Y)}{W_{n}^{\de}(X|X')}\right)\in J \right\}=
-\inf_{q\in J\cap (c_-, c_+)} {\mathcal I}^{\de}(q)
$$
where ${\mathcal I}^{\de}$ is the Legendre transform of $\fee^\de$.
\ep

\br
In the case $\{ Y_t:t\geq 0\}$ is the time reversed process of $\{ X_t:t\geq 0\}$, the cumulant generating function
function $\mathcal{E} (p)$
satisfies the so-called fluctuation theorem symmetry
\[
\mathcal{E} (p) = \mathcal{E} (-1-p)\,.
\]
The large deviation result of Theorem \ref{largedev} then gives  that the entropy production
estimated via waiting times of a discretized version of the process has the same
symmetry in its cumulant generating function for $p\in [0,1]$.
\er

\subsection{Central limit theorem}

\bt\label{clt}
For all $\delta>0$,
$$
\frac1{\sqrt{n}} \left( \log \left(\frac{W_{n}^{\delta}
  (X|Y)}{W_{n}^{\delta} (X|X')}\right)-n\s(\pee|\peet)
\right)
$$
converges in distribution to a normal law
$\mathcal{N}(0,\si_\delta^2)$, where
$$
\si_\delta^2=\lim_{n\to\infty} \frac{1}{n}
\textup{Var}\left(\log\left(\frac{\pee^\delta(X_1^n)}{\peet^\delta(X_1^n)} \right)\right)\,\cdot
$$
Moreover
$$
\lim_{\delta\to 0} \frac1{\delta^2}\si_\delta^2 = \theta^2
$$
where
$$
\theta^2 = \lim_{t\to\infty}\frac{1}{t} \textup{Var}
\left(\log\left(\frac{d\pot}{d\pott} \right)\right)\,\cdot
$$
\et

\bpr
First we claim that for all $\delta >0$
\be\label{l2est}
\lim_{n\to\infty}\frac1n\E_{\pee^\delta\otimes\peet^\delta\otimes\pee^\delta}
\left(\log\left(\frac{W_{n}^{\delta} (X|Y)}{W_{n}^{\delta} (X|X')}\right)
-\sum_{i=1}^n \log\frac{p^X_\delta (X_i,X_{i+1})}{p^Y_\delta (X_i,X_{i+1})}\right)^2
=0\, .
\ee
This follows from the exponential law, as is shown in \cite{cr}, proof
of Theorem 2.

Equation \eqref{l2est} implies that a CLT for $\left(\log W_{n}^{\delta} (X|Y) - \log W_{n}^{\delta} (X|X')\right)$
is equivalent to a CLT for $\sum_{i=1}^n \log\frac{p^X_\delta (X_i,X_{i+1})}{p^Y_\delta (X_i,X_{i+1})}$
and the variances of the asymptotic normals are equal.
For $\delta$ fixed,
$\sum_{i=1}^n \log\frac{p^X_\delta (X_i,X_{i+1})}{p^Y_\delta (X_i,X_{i+1})}$
satisfies the CLT (for $\delta>0$ fixed, $X_i$ is a discrete-time
ergodic Markov chain), so the only thing left is the claimed limiting behavior for the variance, as $\delta\to 0$.

As in the proof of the large deviation theorem, we first develop up to order
$\delta$:
\beq
&&\sum_{i=1}^n \log\frac{p^X_\delta (X_i,X_{i+1})}{p^Y_\delta (X_i,X_{i+1})}
\nonumber\\
&=&
\sum_{i=1}^n \1(X_i=X_{i+1}) \delta (\ct (X_i)-c(X_i)) + \sum_{i=1}^n \1(X_i\not = X_{i+1})
\log\frac{c(X_i)}{\ct(X_i)}
\nonumber\\
&=:&\sum_{i=1}^n (\xi_i^\delta + \zeta_i^\delta)\,.
\nonumber
\eeq
It is then sufficient to verify that
$$
\lim_{\delta\to 0}\lim_{n\to\infty}\frac1{n\delta}\sum_{i=1}^n
\E\left( \left(\xi_i^\delta- \int_{i\delta}^{(i+1)\delta} (\ct(X_s)-c(X_s))ds\right)^2 +
\left(\zeta_i^\delta-\int_{i\delta}^{(i+1)\delta} \log \frac{c(X_s)}{\ct(X_s)} dN_s\right)^2\right)
=0
$$
which is an analogous computation with Poisson random variables as the one used in the proof of
Theorem \ref{largedev}.
\epr

\section{Shadowing a given trajectory}\label{IMI}

Let $\gamma\in D([0,\infty),X)$ be a given trajectory.
The jump process associated to $\gamma$ is defined by
$$
N_t (\gamma) = \sum_{0\leq s\leq t}
\1(\gamma_{s^-}\not=\gamma_{s^+})\, .
$$
For a given $\delta>0$, define the ``jump times" of the $\delta$-discretization of $\gamma$:
$$
\Sigma^\delta_n (\gamma)  = \{ i\in \{ 1,\ldots,n\}:
\gamma_{(i-1)\delta} \not= \gamma_{i\delta}\}\, .
$$
For the Markov process $\{X_t, t\geq 0\}$ with generator
$$
Lf(x) = \sum_{y\in A} c(x)p(x,y) (f(y)-f(x))
$$
define the hitting time
$$
T^\delta_{n} (\gamma|X) =\inf \{k\geq 0 :
(X_{k\delta}, \ldots,X_{(k+n)\delta})= \gamma(\de,\ldots,n\de)\}\, .
$$
In words this is the first time after which the $\delta$-discretization of the process
imitates the $\delta$-discretization of the given trajectory $\gamma$ during $n$ time-steps.
For fixed $\delta >0$, the process $\{ X_{n\delta};n\in\N\}$ is an ergodic
discrete-time Markov chain for which we can apply the results of
\cite{abadi} for hitting times. More precisely there exist
$0<\Lambda_1<\Lambda_2<\infty$ and $C,c,\alpha >0$ such that for all
$\gamma$, $n\in\N$, there exists $\Lambda_1<\lambda_n^\gamma<\Lambda_2$ such that
\be\label{expoo}
\left|\pee \Big( T^\delta_{n} (\gamma|X) \pee(
X_{i\delta}=\gamma_{i\delta} \ ,  \forall i=1,\ldots,n+1)
>t\Big)-e^{-\lambda_n^\gamma t}\right|\leq Ce^{-ct}e^{-\alpha n}\,.
\ee

As a consequence of \eqref{expoo} we have
\bp\label{gammaprop}
For all $\delta>0$,
there exist $\kappa_1,\kappa_2 >0$ such that
for all $\gamma\in D([0,\infty), X)$, $\pee\otimes\peet$ eventually almost surely
$$
-\kappa_1\log n\leq \log \left(T^\delta_n (\gamma|Y) \peet
(Y_\delta=\gamma_\delta,\ldots,Y_{n\delta}=\gamma_{n\delta})\right)\leq
\log(\log n^{\kappa_2})
$$
and
$$
-\kappa_1\log n\leq \log \left(T^\delta_n (\gamma|X) \pee
(X_\delta=\gamma_{\delta},\ldots,X_{n\delta}=\gamma_{n\delta})\right)\leq
\log(\log n^{\kappa_2})\,.
$$
\ep

Therefore, for $\delta>0$ fixed, we arrive at
\beq
\label{shit}
\lefteqn{\log T^\delta_{n} (\gamma|X)= }\\
\sum_{i\in\Sigma^\delta_n (\gamma)}
&& {} \log (\de c(\gamma_{(i-1)\de}) p(\gamma_{(i-1)\de},\gamma_{i\de}))
+
\sum_{i\in \{ 1,\ldots,n\}\setminus \Sigma^\delta_n (\gamma)}
\log(1-\de c(\gamma_{(i-1)\de})) + o(n)\,.
\nonumber
\eeq
The presence of the $\log (\de)$ term in the rhs of \eqref{shit} causes the same problem
as we have encountered in Section \ref{naive}. Therefore, we have to subtract another quantity such that
the $\log (\de)$ term is canceled.
In the spirit of what we did with the waiting times, we subtract $\log T_n^\de (\gamma|Y)$, where
$\{ Y_t:t\geq 0\}$ is another independent Markov process with generator
$$
Lf(x) = \sum_{y\in A} \tilde{c}(x)\tilde{p}(x,y) (f(y)-f(x))\,.
$$
We then arrive at
\be\label{quot}
\log \frac{T^\delta_{n} (\gamma|X)}{ T^\delta_{n} (\gamma|X)}
=
\sum_{i\in\Sigma^\delta_n (\gamma)}\log \frac{c(\gamma_{(i-1)\de}) p(\gamma_{(i-1)\de},\gamma_{i\de})}{\ct (\gamma_{(i-1)\de}) \pt
(\gamma_{(i-1)\de},\gamma_{i\de})}
+
\sum_{i\in \{ 1,\ldots,n\}\setminus \Sigma^\delta_n (\gamma)}
\log\frac{(1-\de c(\gamma_{(i-1)\de}))}{(1-\de \ct (\gamma_{(i-1)\de}))} + o(n)
\ee

We then have the following law of large numbers
\bt Let $\pee$ (resp.\ $\peet$) denote the stationary path space measure of $\{ X_t:t\geq 0\}$
(resp. $\{ Y_t:t\geq 0\}$ and
let $\gamma\in D([0,\infty), X)$ be a fixed trajectory. We then have
$\pee\otimes\peet$-almost surely:
\beq\label{quotin}
& &\lim_{\de\to 0} \lim_{n\to\infty}
\frac{1}{n\de}\left(
\log\frac{ T_n^{\de} (\gamma|Y)}{ T^{\de}_n (\gamma|X)} -
\int_0^{n\de}
\log\frac{c(\gamma_s) p(\gamma_{s^-},\gamma_{s^+})}{\tilde{c}(\gamma_s) \pt (\gamma_{s^-},\gamma_{s^+})}
dN_s (\gamma)
-\int_0^{n\de}
(\ct (\gamma_s) -c(\gamma_s))ds
\right)
\nonumber\\
&&= 0\,.
\eeq
Moreover, if $\gamma$ is chosen according to a stationary ergodic measure $\mathbb{Q}$
on path-space, then $\mathbb{Q}$-almost surely
\be
\lim_{\de\to 0}\lim_{n\to\infty}\frac{1}{n\de}\left(
\log T_n^{\de} (\gamma|Y)-\log T^{\de}_n (\gamma|X) \right)
=
\sum_{x,y\in A} q(x,y) \log \frac{c(x) p(x,y)}{\ct (x) \pt (x,y)}
+
\sum_{x\in A} q(x)(\ct (x)-c(x))
\ee
where
\beq
q(x,y) &=& \lim_{t\to\infty}\E_{\qak}\left( \frac{N^{xy}_t}{t}\right)
\nonumber\\
q(x) &=& \qak (\gamma_0 = x)
\nonumber
\eeq
and where $N^{xy}_t (\gamma)$ denotes the number of jumps from $x$ to $y$ of the trajectory $\gamma$
in the time-interval
$[0,t]$.
\et
\bpr
Using proposition \ref{gammaprop}, we use the same proof as
that of Theorem \ref{lln}, and use
that the sums in the rhs of \eqref{quot} is
up to order $\delta^2$  equal to the integrals appearing in the lhs
of \eqref{quotin}.
The other assertions of the theorem follow from the ergodic theorem.
\epr
\br
If we choose $\gamma$ according to the path space measure $\pee$, i.e.,
$\gamma$ is a ``typical" trajectory of the process $\{ X_t:t \geq 0\}$, and choose
$p(x,y)=\pt (x,y)$, then we recover the limit of the law of large numbers
for waiting times (Theorem \ref{lln}):
\[
\lim_{\de\to 0}\lim_{n\to\infty}\frac{1}{n\de}\left(
\log T_n^{\de} (\gamma|Y)-\log T^{\de}_n (\gamma|X) \right)
=
\sum_{x}\mu (x)c(x) \log\frac{c(x)}{\ct (x)} +
\sum_{x}\mu(x) (\ct(x)-c(x))= s(\pee|\peet)
\]

\er


\begin{thebibliography}{20}

\bibitem{abadi}
M. Abadi, {\em Exponential approximation for hitting times in mixing
processes}, Math. Phys. Electron. J. {\bf 7} (2001).

\bibitem{acrv}
M. Abadi, J.-R. Chazottes F. Redig and E. Verbitskiy,
{\em Exponential distribution for the occurrence
of rare patterns in Gibbsian random fields},
Commun. Math. Phys. {\bf 246} no. 2 (2004), 269--294.

\bibitem{cr}
J.-R.\ Chazottes and F.\ Redig,
{\em Testing the irreversibility of a Gibbsian process via hitting and
  return times}, Nonlinearity (2005), {\bf 18}, 2477--2489.

\bibitem{DZ}
A. Dembo and O. Zeitouni, Large deviation techniques and applications,
Springer, (1998).

\bibitem{GS}
I.I. Gihman, A.V. Skorohod,
The theory of stochastic processes. II.
Die Grundlehren der Mathematischen Wissenschaften {\bf 218},
Springer, 1975.

\bibitem{JQQ}
D.-Q. Jiang, M. Qian, M.-P. Qian,
{\em Mathematical theory of nonequilibrium steady states.
On the frontier of probability and dynamical systems}.
Lecture Notes in Mathematics {\bf 1833}, Springer, 2004.

\bibitem{maes}
C. Maes, {\em The fluctuation theorem as a Gibbs property},
J. Stat. Phys. {\bf 95}, 367-392, (1999).

\bibitem{PS}
D. Plachky, J. Steinebach,
{\em A theorem about probabilities of large deviations with an application to queuing theory},
Period. Math. Hungar. {\bf 6} (1975), no. 4, 343--345.

\bibitem{shields}
P.C. Shields,
The ergodic theory of discrete sample paths.
Graduate Studies in Mathematics {\bf 13},
American Mathematical Society, Providence, RI, 1996.

\bibitem{stroock}
D. Stroock,
An introduction to Markov processes.
Graduate Texts in Mathematics {\bf 230},
Springer, 2005.

\bibitem{varadhan}
S.R.S. Varadhan,
Large deviations and applications. Philadelphia: Society for Industrial and Applied
Mathematics, 1984.

\end{thebibliography}
\end{document}